\documentclass[aos,MSNbibl,nameyear,dvips]{arximspdf}

\doi{10.1214/12-AOS985} 
\referstodoi{10.1214/11-AOS949}
\volume{40}
\issue{4}
\pubyear{2012}
\firstpage{1989}
\lastpage{1996}

\makeatletter
\newtheorem{theorem}{Theorem}
\newproclaim{remark}{Remark}
\makeatother

\begin{document}
\begin{frontmatter}

\title{Discussion: Latent variable graphical model selection via
convex optimization\thanksref{T1}}
\thankstext{T1}{Supported in part by NSF Career Award DMS-06-45676 and
NSF FRG Grant DMS-08-54975.}
\runtitle{Comment}

\begin{aug}
\author[A]{\fnms{Zhao} \snm{Ren}\ead[label=e1]{zhao.ren@yale.edu}}
\and
\author[A]{\fnms{Harrison H.} \snm{Zhou}\corref{}\ead[label=e2]{huibin.zhou@yale.edu}}
\runauthor{Z. Ren and H. H. Zhou}
\affiliation{Yale University}
\address[A]{Department of Statistics\\
Yale University\\
New Haven, Connecticut 06511\\
USA\\
\printead{e1}\\
\phantom{E-mail:\ }\printead*{e2}} 
\end{aug}

\received{\smonth{2} \syear{2012}}



\end{frontmatter}
%

\section{Introduction}

We would like to congratulate the authors for their refreshing contribution
to this high-dimensional latent variables graphical model selection problem.
The problem of covariance and concentration matrices is fundamentally
important in several classical statistical methodologies and many
applications. Recently, sparse concentration matrices estimation has
received considerable attention, partly due to its connection to sparse
structure learning for Gaussian graphical models. See, for example,
Meinshausen and B\"{u}hlmann (\citeyear{MeiBuh06}) and Ravikumar et al. (\citeyear{Ravetal11}). Cai,
Liu and Zhou (\citeyear{CaiLiuZho}) considered rate-optimal estimation.

The authors extended the current scope to include latent variables. They
assume that the fully observed Gaussian graphical model has a naturally
sparse dependence graph. However, there are only partial observations
available for which the graph is usually no longer sparse. Let $X$ be $%
( p+r)$-variate Gaussian with a sparse concentration matrix $%
S_{( O,H) }^{\ast}$. We only\vspace*{1pt} observe $X_{O}$, $p$ out of the
whole $p+r$ variables, and denote its covariance matrix by $\Sigma
_{O}^{\ast}$. In this case, usually the $p\times p$ concentration
matrix $%
( \Sigma_{O}^{\ast}) ^{-1}$ are not sparse. Let $S^{\ast}$ be
the concentration matrix of observed variables conditioned on latent
variables, which is a submatrix\vspace*{2pt} of $S_{( O,H) }^{\ast}$ and
hence has a sparse structure, and let $L^{\ast}$ be the summary of the
marginalization over the latent variables and its rank corresponds to the
number of latent variables $r$ for which we usually assume it is small. The
authors observed $( \Sigma_{O}^{\ast}) ^{-1}$ can be decomposed
as the difference of the sparse matrix $S^{\ast}$ and the rank $r$
matrix $%
L^{\ast}$, that is, $( \Sigma_{O}^{\ast}) ^{-1}=S^{\ast}-L^{\ast
}$. Then following traditional wisdoms, the authors naturally proposed a
\textit{regularized maximum likelihood approach} to estimate both the sparse
structure $S^{\ast}$ and the low-rank part $L^{\ast}$,
\[
\min_{( S,L) :S-L\succ0,L\succeq0}\operatorname{tr}\bigl(
( S-L) \Sigma_{O}^{n}\bigr) -\log\det( S-L) +\chi
_{n}\bigl( \gamma\Vert S\Vert_{1}+\operatorname{tr}( L)
\bigr),
\]
where $\Sigma_{O}^{n}$ is the sample covariance matrix, $\Vert
S\Vert_{1}=\sum_{i,j}\vert s_{ij}\vert$, and $\gamma$
and $\chi_{n}$ are regularization tuning parameters. Here $\operatorname{tr}
( L) $ is the trace of $L$. The notation $A\succ0$ means $A$ is
positive definite, and $A\succeq0$ denotes that $A$ is nonnegative.

There is an obvious identifiability problem if we want to estimate both the
sparse and low-rank components. A matrix can be both sparse and low
rank. By
exploring the geometric properties of the tangent spaces for sparse and
low-rank components, the authors gave a beautiful sufficient condition for
identifiability, and then provided very much involved theoretical
justifications based on the sufficient condition, which is beyond our
ability to digest them in a short period of time in the sense that we don't
fully understand why those technical assumptions were needed in the analysis
of their approach. Thus, we decided to look at a relatively simple but
potentially practical model, with the hope to still capture the essence of
the problem, and see how well their regularized procedure works. Let $%
\Vert\cdot\Vert_{1\rightarrow1}$ denote the matrix $l_{1}$
norm, that is, $\Vert S\Vert_{1\rightarrow1}=\max_{1\leq i\leq
p}\sum_{j=1}^{p}\vert s_{ij}\vert$. We assume that $S^{\ast}$
is in the following uniformity class:
\begin{eqnarray}\label{Sstar}
\mathcal{U}( s_{0}( p) ,M_{p}) &=&\Biggl\{ S=(
s_{ij}) \dvtx S\succ0,\Vert S\Vert_{1\rightarrow1}\leq M_{p},\nonumber\\[-8pt]\\ [-8pt]
&&\hspace*{17pt}\hphantom{\Biggl\{}{}\max_{1\leq i\leq p}\sum_{j=1}^{p}\mathbf{1}\{ s_{ij}\neq
0\} \leq s_{0}( p) \Biggr\} ,\nonumber
\end{eqnarray}
where we allow $s_{0}( p) $ and $M_{p}$ to grow as $p$ and $n$
increase. This uniformity class was considered in Ravikumar et al. (\citeyear{Ravetal11})
and Cai, Liu and Luo (\citeyear{CaiLiuLuo11}). For the low-rank matrix $L^{\ast}$, we assume
that the effect of marginalization over the latent variables spreads out,
that is, the low-rank matrix $L^{\ast}$ has row/column spaces that are not
closely aligned with the coordinate axes to resolve the identifiability
problem. Let the eigen-decomposition of $L^{\ast}$ be as follows:%
%
\begin{equation}\label{Lstar}
L^{\ast}=\sum_{i=1}^{r_{0}( p) }\lambda_{i}u_{i}u_{i}^{T},
\end{equation}
where $r_{0}(p)$ is the rank of $L^{\ast}$. We assume that there
exists a
universal constant $c_{0}$ such that $\Vert u_{i}\Vert_{\infty
}\leq\sqrt{\frac{c_{0}}{p}}$ for all $i$, and $\Vert L^{\ast
}\Vert_{1\rightarrow1}$ is bounded by $M_{p}$ which can be shown to
be bounded by $c_{0}r_{0}$. A similar incoherence assumption on $u_{i}$ was
used in Cand\`{e}s and Recht (\citeyear{CanRec09}). We further assume that
%
\begin{equation}\label{eigenbd}
\lambda_{\max}(\Sigma_{O}^{\ast})\leq M\quad \mbox{and}\quad\lambda_{\min
}(\Sigma_{O}^{\ast})\geq1/M
\end{equation}
for some universal constant $M$.

As discussed in the paper, the goals in latent variable model selection are
to obtain the sign consistency for the sparse matrix $S^{\ast}$ as
well\vadjust{\goodbreak} as
the rank consistency for the low-rank semi-positive definite
matrix
$L^{\ast
}$. Denote the minimum magnitude of nonzero entries of $S^{\ast}$ by $%
\theta$, that is, $\theta=\min_{i,j}\vert s_{ij}\vert\mathbf{1}%
\{ s_{ij}\neq0\} $, and the minimum nonzero eigenvalue of $%
L^{\ast}$ by $\sigma$, that is, $\sigma=\min_{1\leq i\leq
r_{0}}\lambda_{i}$%
. To obtain theoretical guarantees of consistency results for the model
described in (\ref{Sstar}), (\ref{Lstar}) and (\ref{eigenbd}), in addition
to the strong irrepresentability condition which seems to be difficult to
check in practice, the authors require the following assumptions (by a
translation of the conditions in the paper to this model) for $\theta
,\sigma$ and $n$:

\begin{longlist}
\item[(1)] $\theta\gtrsim\sqrt{p/n},$ which is needed even when $s_{0}(p)$
is constant;

\item[(2)] $\sigma$ $\gtrsim$ $s_{0}^{3}( p) \sqrt{{p/n}}$
under the additional strong assumptions on the Fisher information
matrix $%
\Sigma_{O}^{\ast}\otimes\Sigma_{O}^{\ast}$ (see the footnote for
Corollary 4.2);

\item[(3)] $n\gtrsim s_{0}^{4}( p) p$.
\end{longlist}

However, for sparse graphical model selection without latent variables,
either the $l_{1}$-regularized maximum likelihood approach [see
Ravikumar et al.
(\citeyear{Ravetal11})] or CLIME [see Cai, Liu and Luo (\citeyear{CaiLiuLuo11})] can be shown to be sign
consistent if the minimum magnitude nonzero entry of concentration
matrix $%
\theta$ is at the order of $\sqrt{( \log p) /n}$ when $M_{p}$ is
bounded, which inspires us to study rate-optimalites for this latent
variables graphical model selection problem. In this discussion, we propose
a procedure to obtain an algebraically consistent estimate of the latent
variable Gaussian graphical model under a much weaker condition on both
$%
\theta$ and~$\sigma$. For example, for a wide range of $s_{0}(
p) $, we only require $\theta$ is at the order of $\sqrt{( \log
p) /n}$ and $\sigma$ is at the order of $\sqrt{p/n}$ to consistently
estimate the support of $S^{\ast}$ and the rank of $L^{\ast}$. That means
the \textit{regularized maximum likelihood approach} could be far from being
optimal, but we don't know yet whether the suboptimality is due to the
procedure or their theoretical analysis.

\section{Latent variable model selection consistency}

In this section we propose a procedure to obtain an algebraically
consistent estimate of the latent variable Gaussian graphical model. The
condition on $\theta$ to recover the support of $S^{\ast}$ is reduced to
that in Cai, Liu and Luo (\citeyear{CaiLiuLuo11}) which studied sparse graphical model
selection without latent variables, and the condition on $\sigma$ is just
at an order of $\sqrt{p/n}$, which is smaller than $s_{0}^{3}( p)
\sqrt{{p/n}}$ assumed in the paper when $s_{0}( p) \rightarrow
\infty$. When $M_{p}$ is bounded, our results can be shown to be
rate-optimal by lower bounds stated in Remarks \ref{lw1} and \ref
{lw2} for
which we are not giving proofs due to the limitation of the space.

\subsection{\texorpdfstring{Sign consistency procedure of $S^{\ast}$}
{Sign consistency procedure of S*}}

We propose a CLIME-like estimator of $S^{\ast}$ by solving the following
linear optimization problem:
\[
\min\Vert S\Vert_{1}\qquad\mbox{subject to}\qquad\Vert\Sigma
_{O}^{n}S-I\Vert_{\infty}\leq\tau_{n},\qquad S\in\mathbb{R}^{p\times p},
\]
where $\Sigma_{O}^{n}=( \widetilde{\sigma}_{ij}) $ is the
sample covariance matrix. The tuning parameter $\tau_{n}$ is chosen as
$%
\tau_{n}=C_{1}M_{p}\sqrt{\frac{\log p}{n}}$\vadjust{\goodbreak} for some large constant
$C_{1}$%
. Let $\hat{S}_{1}=( \hat{s}_{ij}^{1}) $ be the solution. The
CLIME-like estimator $\hat{S}=( \hat{s}_{ij}) $ is obtained by
symmetrizing $\hat{S}_{1}$ as follows:
\[
\hat{s}_{ij}=\hat{s}_{ji}=\hat{s}_{ij}^{1}\mathbf{1}\{ \vert\hat{s
}_{ij}^{1}\vert\leq\hat{s}_{ji}^{1}\} +\hat{s}%
_{ji}^{1}\mathbf{1}\{ \vert\hat{s}_{ij}^{1}\vert>\hat{s}%
_{ji}^{1}\}.
\]
In other words,\vspace*{-3pt} we take the one with smaller magnitude between $\hat{s}
_{ij}^{1}$ and $\hat{s}_{ji}^{1}$. We define a thresholding estimator $
\tilde{S}=( \tilde{s}_{ij}) $ with
%
\begin{equation}\label{St}
\tilde{s}_{ij}=\tilde{s}_{ij}\mathbf{1}\{ \vert\tilde{s}_{ij}\vert
>9M_{p}\tau_{n}\}
\end{equation}
to estimate the support of $S^{\ast}$.

\begin{theorem}\label{Support Recovery}
Suppose that $S^{\ast}\in\mathcal{U}(
s_{0}( p) ,M_{p}) $,
%
\begin{equation}\label{assump1}
\sqrt{{(\log p)/n}}=o(1)\quad\mbox{and}\quad\Vert L^{\ast}\Vert
_{\infty}\leq M_{p}\tau_{n}.
\end{equation}
With probability greater than $1-C_{s}p^{-6}$ for some constant $C_{s}$
depending on $M$ only, we have%
\[
\Vert\hat{S}-S^{\ast}\Vert_{\infty}\leq9M_{p}\tau_{n}.
\]
Hence, if the minimum magnitude of nonzero entries $\theta>18M_{p}\tau_{n},$
we obtain the sign consistency $\operatorname{sign}( \tilde{S}) =\operatorname{sign}(
S^{\ast}) $. In particular, if $M_{p}$ is in the constant level, then
to consistently recover the support of $S^{\ast}$, we only need that $%
\theta\asymp\sqrt{(\log p)/n}.$
\end{theorem}

\begin{pf}
The proof is similar to Theorem $7$ in Cai, Liu and Luo (\citeyear{CaiLiuLuo11}). The
sub-Gaussian condition with spectral norm upper bound $M$ implies that each
empirical covariance $\widetilde{\sigma}_{ij}$ satisfies the following
large deviation result:%
\[
\mathbb{P}( \vert\widetilde{\sigma}_{ij}-\sigma_{ij}\vert
>t) \leq C_{s}\exp\biggl( -\frac{8}{C_{2}^{2}}nt^{2}\biggr)\qquad\mbox{for }\vert t\vert\leq\phi,
\]
where $C_{s},C_{2}$ and $\phi$ only depend on $M$. See, for example,
Bickel and Levina (\citeyear{BicLev08}). In particular, for $t=C_{2}\sqrt{( \log
p) /n}$ which is less than $\phi$ by our assumption, we have
%
\begin{equation}\label{TailPr1}
\mathbb{P}( \Vert\Sigma_{O}^{\ast}-\Sigma_{O}^{n}\Vert
_{\infty}>t) \leq\sum_{i,j}\mathbb{P}( \vert\widetilde{%
\sigma}_{ij}-\sigma_{ij}\vert>t) \leq p^{2}\cdot C_{s}p^{-8}.
\end{equation}

Let
\[
A=\bigl\{ \Vert\Sigma_{O}^{\ast}-\Sigma_{O}^{n}\Vert_{\infty
}\leq C_{2}\sqrt{{(\log p)/n}}\bigr\}.
\]
Equation (\ref{TailPr1})\ implies $\mathbb{P}( A) \geq
1-C_{s}p^{-6}$. On event $A$, we will show
%
\begin{equation} \label{med}
\Vert( S^{\ast}-L^{\ast}) -\hat{S}_{1}\Vert
_{\infty}\leq8M_{p}\tau_{n},
\end{equation}
which immediately yields%
\[
\Vert S^{\ast}-\hat{S}\Vert_{\infty}\leq\Vert(
S^{\ast}-L^{\ast}) -\hat{S}_{1}\Vert_{\infty}+\Vert
L^{\ast}\Vert_{\infty}
\leq8M_{p}\tau_{n}+M_{p}\tau
_{n}=9M_{p}\tau_{n}.
\]

Now we establish equation (\ref{med}). On event $A$, for some large constant
$C_{1}\geq2C_{2}$, the choice of $\tau_{n}$ yields
%
\begin{equation}\label{SparseTuning}
2M_{p}\Vert\Sigma_{O}^{\ast}-\Sigma_{O}^{n}\Vert_{\infty
}\leq\tau_{n}.
\end{equation}
By the matrix $l_{1}$ norm assumption, we could obtain that
%
\begin{equation}\label{SandSigma}
\Vert( \Sigma_{O}^{\ast}) ^{-1}\Vert_{1\rightarrow
1}\leq\Vert S^{\ast}\Vert_{1\rightarrow1}+\Vert L^{\ast
}\Vert_{1\rightarrow1}
\leq2M_{p}.
\end{equation}
From (\ref{SparseTuning}) and (\ref{SandSigma}) we have%
\begin{eqnarray*}
\Vert\Sigma_{O}^{n}( S^{\ast}-L^{\ast}) -I\Vert
_{\infty}&=&\Vert( \Sigma_{O}^{n}-\Sigma_{O}^{\ast})
( \Sigma_{O}^{\ast}) ^{-1}\Vert_{\infty}\\
&\leq&\Vert
\Sigma_{O}^{n}-\Sigma_{O}^{\ast}\Vert_{\infty}\Vert(
\Sigma_{O}^{\ast}) ^{-1}\Vert_{1\rightarrow1}\leq\tau_{n},
\end{eqnarray*}
which implies%
\begin{eqnarray}\label{condition1}
&&\Vert\Sigma_{O}^{n}( S^{\ast}-L^{\ast}) -\Sigma_{O}^{n}%
\hat{S}_{1}\Vert_{\infty}\nonumber\\ [-8pt]\\ [-8pt]
&&\qquad\leq\Vert\Sigma_{O}^{n}(
S^{\ast}-L^{\ast}) -I\Vert_{\infty}+\Vert\Sigma
_{O}^{n}\hat{S}_{1}-I\Vert_{\infty}
\leq2\tau_{n}.\nonumber
\end{eqnarray}
From the definition of $\hat{S}_{1}$ we obtain that%
%
\begin{equation}\label{L1normCom1}
\Vert\hat{S}_{1}\Vert_{1\rightarrow1}\leq\Vert S^{\ast
}-L^{\ast}\Vert_{1\rightarrow1}\leq2M_{p},
\end{equation}
which, together with equations (\ref{SparseTuning}) and (\ref{condition1}),
implies%
\begin{eqnarray*}
&&\bigl\Vert\Sigma_{O}^{\ast}\bigl( ( S^{\ast}-L^{\ast}) -\hat{%
S}_{1}\bigr) \bigr\Vert_{\infty}\\
 &&\quad\leq\Vert\Sigma_{O}^{n}(
S^{\ast}-L^{\ast}) -\hat{S}_{1}\Vert_{\infty}+\bigl\Vert
( \Sigma_{O}^{\ast}-\Sigma_{O}^{n}) \bigl( ( S^{\ast
}-L^{\ast}) -\hat{S}_{1}\bigr) \bigr\Vert_{\infty} \\
&&\quad\leq2\tau_{n}+\Vert\Sigma_{O}^{n}-\Sigma_{O}^{\ast}\Vert
_{\infty}\Vert( S^{\ast}-L^{\ast}) -\hat{S}%
_{1}\Vert_{1\rightarrow1} \\
&&\quad\leq2\tau_{n}+4M_{p}\Vert\Sigma_{O}^{n}-\Sigma_{O}^{\ast
}\Vert_{\infty}\leq4\tau_{n}.
\end{eqnarray*}
Thus, we have%
\[
\Vert( S^{\ast}-L^{\ast}) -\hat{S}_{1}\Vert
_{\infty}\leq\Vert( \Sigma_{O}^{\ast}) ^{-1}\Vert
_{1\rightarrow1}\bigl\Vert\Sigma_{O}^{\ast}\bigl( ( S^{\ast
}-L^{\ast}) -\hat{S}_{1}\bigr) \bigr\Vert_{\infty}\\
\leq8M_{p}\tau
_{n}.
\]
\upqed\end{pf}

\begin{remark}
By the choice of our $\tau_{n}$ and the eigen-decomposition of
$L^{\ast}$,
the condition $\Vert L^{\ast}\Vert_{\infty}\leq M_{p}\tau_{n}$
holds when $r_{0}(p)C_{0}/p\leq C_{1}M_{p}^{2}\sqrt{( \log p) /n}$%
, that is, $p^{2}\log p\gtrsim nr_{0}^{2}(p)M_{p}^{-4}$. If $M_{p}$ is slowly
increasing (e.g., $p^{1/4-\tau}$ for any small \mbox{$\tau>0$}), the
minimum requirement $\theta\asymp M_{p}^{2}\sqrt{( \log p) /n}$
is weaker than $\theta\gtrsim\sqrt{p/n}$ required in Corollary $4.2$.
Furthermore, it can be shown that the optimal rate of minimum magnitude of
nonzero entries for sign consistency is $\theta\asymp M_{p}\sqrt
{(\log p)/n}
$ as in Cai, Liu and Zhou (\citeyear{CaiLiuZho}).
\end{remark}

\begin{remark}\label{lw1}
Cai, Liu and Zhou (\citeyear{CaiLiuZho}) showed the minimum requirement for
$%
\theta$, $\theta\asymp M_{p}\sqrt{(\log p)/n}$ is necessary for sign
consistency for sparse concentration matrices. Let $\mathcal{U}_{S}(c)$
denote the class of concentration matrices defined in (\ref{Sstar})
and (\ref{Lstar}), satisfying assumption (\ref{assump1}) and $\theta
>cM_{p}\sqrt{%
(\log p)/n}$. We can show that there exists some constant $c_{1}>0$ such
that for all $0<c<c_{1},$%
\[
\lim_{n\rightarrow\infty}\inf_{( \hat{S},\hat{L}) }\sup_{%
\mathcal{U}_{S}(c)}\mathbb{P}\bigl( \operatorname{sign}( \hat{S}) \neq
\operatorname{sign}( S^{\ast}) \bigr) >0,
\]
similar to Cai, Liu and Zhou (\citeyear{CaiLiuZho}).
\end{remark}

\subsection{\texorpdfstring{Rank Consistency Procedure of $L^{\ast}$}
{Rank Consistency Procedure of L*}}

In this section we propose a procedure to estimate $L^{\ast}$ and its rank.
We note that with high probability $\Sigma_{O}^{n}$ is invertible, then
define $\hat{L}=( \Sigma_{O}^{n}) ^{-1}-\tilde{S},$ where $%
\tilde{S}$ is defined in (\ref{St}). Denote the eigen-decomposition
of $\hat{%
L}$ by $\sum_{i=1}^{p}\lambda_{i}(\hat{L})\upsilon_{i}\upsilon_{i}^{T}$,
and let $\lambda_{i}(\tilde{L})=\lambda_{i}(\hat{L})1\{ \lambda
_{i}(%
\hat{L})>C_{3}\sqrt{\frac{p}{n}}\} $, where constant $C_{3}$ will be
specified later. Define $\tilde{L}=\sum_{i=1}^{p}\lambda_{i}(\tilde
{L}%
)\upsilon_{i}\upsilon_{i}^{T}$. The following theorem shows that estimator
$\tilde{L}$ is a consistent estimator of $L^{\ast}$ under the
spectral norm
and with high probability $\operatorname{rank}( L^{\ast}) =\operatorname{rank}( \tilde{L}%
) $.

\begin{theorem}\label{Rank Consistancy}
Under the conditions in Theorem \ref{Support
Recovery}, we assume that
%
\begin{equation}\label{assmp}
\sqrt{\frac{p}{n}}\leq\frac{1}{16\sqrt{2}M^{2}} \quad\mbox{and}\quad%
M_{p}^{2}s_{0}(p)\leq\sqrt{\frac{p}{\log p}}.
\end{equation}
Then there exists some constant $C_{3}$ such that
\[
\Vert\hat{L}-L^{\ast}\Vert\leq C_{3}\sqrt{\frac{p}{n}}
\]
with probability greater than $1-2e^{-p}-C_{s}p^{-6}$. Hence, if
$\sigma
>2C_{3}\sqrt{\frac{p}{n}},$ we have $\operatorname{rank}( L^{\ast}) =\operatorname{rank}(
\tilde{L}) $ with high probability.
\end{theorem}

\begin{pf}
From Corollary $5.5$ of the paper and our assumption on the sample size,
we have%
\[
\mathbb{P}\biggl( \Vert\Sigma_{O}^{\ast}-\Sigma_{O}^{n}\Vert
\geq\sqrt{128}M\sqrt{\frac{p}{n}}\biggr) \leq2\exp( -p).
\]
Note that $\lambda_{\min}( \Sigma_{O}^{\ast}) $ $\geq$ $1/M$%
, and $\sqrt{128}M\sqrt{\frac{p}{n}}\leq1/( 2M) $ under the
assumption (\ref{assmp}), then $\lambda_{\min}(\Sigma_{O}^{n})\geq
1/( 2M) $ with high probability, which yields the same rate of
convergence for the concentration matrix, since%
\begin{eqnarray}\label{SpectralBoundInverse}
\Vert( \Sigma_{O}^{\ast}) ^{-1}-( \Sigma
_{O}^{n}) ^{-1}\Vert&\leq&\Vert( \Sigma_{O}^{\ast
}) ^{-1}\Vert\Vert( \Sigma_{O}^{n})
^{-1}\Vert\Vert\Sigma_{O}^{\ast}-\Sigma_{O}^{n}\Vert\nonumber\\ [-8pt]\\ [-8pt]
&\leq&2M^{2}\sqrt{128}M\sqrt{\frac{p}{n}}=16\sqrt{2}M^{3}\sqrt
{\frac{p}{n}}.\nonumber
\end{eqnarray}
From Theorem \ref{Support Recovery} we know%
\[
\operatorname{sign}( \tilde{S}) =\operatorname{sign}( S^{\ast})\quad \mbox{and}\quad
\Vert\tilde{S}-S^{\ast}\Vert_{\infty}\leq9M_{p}\tau_{n}\vadjust{\goodbreak}
\]
with probability greater than $1-C_{s}p^{-6}$. Since $\Vert
B\Vert\leq\Vert B\Vert_{1\rightarrow1}$ for any
symmetric matrix $B$, we then have%
%
\begin{equation}\label{SpectralBoundSparse}
\Vert\tilde{S}-S^{\ast}\Vert\leq\Vert\tilde{S}-S^{\ast
}\Vert_{1\rightarrow1}\leq s_{0}( p) 9M_{p}\tau
_{n}=9C_{1}M_{p}^{2}s_{0}( p) \sqrt{\frac{\log p}{n}}.
\end{equation}
Equations (\ref{SpectralBoundInverse}) and (\ref{SpectralBoundSparse}),
together with the assumption $M_{p}^{2}s_{0}(p)\leq\sqrt{\frac
{p}{\log p}}$,
imply%
\begin{eqnarray*}
\Vert\hat{L}-L^{\ast}\Vert&\leq&\Vert( \Sigma
_{O}^{\ast}) ^{-1}-( \Sigma_{O}^{n}) ^{-1}\Vert
+\Vert\tilde{S}-S^{\ast}\Vert\\
&\leq&16\sqrt{2}M^{3}\sqrt{\frac{p%
}{n}}+9C_{1}M_{p}^{2}s_{0}( p) \sqrt{\frac{\log p}{n}}\leq C_{3}%
\sqrt{\frac{p}{n}}
\end{eqnarray*}
with probability greater than $1-2e^{-p}-C_{s}p^{-6}$.
\end{pf}

\begin{remark}
We should emphasize the fact that in order to consistently estimate the rank
of $L^{\ast}$ we need only\vspace*{-3pt} that $\sigma>2C_{3}\sqrt{\frac{p}{n}}$, which
is smaller than $s_{0}^{3}( p) \sqrt{\frac{p}{n}}$ required in
the paper (see the footnote for Corollary 4.2), as long as $%
M_{p}^{2}s_{0}(p)\leq\sqrt{\frac{p}{\log p}}$. In particular, we don't
explicitly constrain the rank $r_{0}(p)$. One special case is that $M_{p}$
is constant and $s_{0}(p)\asymp p^{1/2-\tau}$ for some small $\tau
>0$, for
which our requirement is $\sqrt{\frac{p}{n}}$ but the assumption in the
paper is at an order of $p^{3( 1/2-\tau) }\sqrt{\frac{p}{n}}$.
\end{remark}

\begin{remark}\label{lw2}
Let $\mathcal{U}_{L}(c)$ denote the class of concentration
matrices defined in (\ref{Sstar}), (\ref{Lstar}) and (\ref{eigenbd}),
satisfying assumptions (\ref{assmp}), (\ref{assump1}) and $\sigma
>c\sqrt{%
\frac{p}{n}}$. We can show that there exists some constant $c_{2}>0$ such
that for all $0<c<c_{2}$,%
\[
\lim_{n\rightarrow\infty}\inf_{( \hat{S},\hat{L}) }\sup_{%
\mathcal{U}_{L}(c)}\mathbb{P}\bigl( \operatorname{rank}( \hat{L}) \neq
\operatorname{rank}( L^{\ast}) \bigr) >0.
\]
The proof of this lower bound is based on a modification of a lower bound
argument in a personal communication of T. Tony Cai (\citeyear{Cai}).
\end{remark}

\section{Concluding remarks and further questions}

In this discussion we attempt to understand optimalities of results in the
present paper by studying a relatively simple model. Our preliminary
analysis seems to indicate that their results in this paper are suboptimal.
In particular, we tend to conclude that assumptions on $\theta$ and
$\sigma$
in the paper can be potentially very much weakened. However, it is not clear
to us whether the suboptimality is due to the methodology or just its
theoretical analysis. We want to emphasize that the preliminary results in
this discussion can be strengthened, but for the purpose of simplicity of
the discussion we choose to present weaker but simpler results to hopefully
shed some light on understanding optimalities in estimation.


%

\printaddresses


\begin{thebibliography}{8}

\bibitem[\protect\citeauthoryear{Bickel and Levina}{2008}]{BicLev08}
\begin{barticle}[mr]
\bauthor{\bsnm{Bickel},~\bfnm{Peter~J.}\binits{P.~J.}} \AND
  \bauthor{\bsnm{Levina},~\bfnm{Elizaveta}\binits{E.}}
(\byear{2008}).
\btitle{Regularized estimation of large covariance matrices}.
\bjournal{Ann. Statist.}
\bvolume{36}
\bpages{199--227}.
\bid{doi={10.1214/009053607000000758}, issn={0090-5364}, mr={2387969}}
\bptok{imsref}%
\end{barticle}
\endbibitem

\bibitem[\protect\citeauthoryear{Cai}{2011}]{Cai}
\begin{bmisc}[auto:STB|2012/05/30|10:51:56]
\bauthor{\bsnm{Cai},~\bfnm{T.~T.}\binits{T.~T.}}
(\byear{2011}).
\bhowpublished{Personal communication}.
\bptok{imsref}%
\end{bmisc}
\endbibitem

\bibitem[\protect\citeauthoryear{Cai, Liu and Luo}{2011}]{CaiLiuLuo11}
\begin{barticle}[mr]
\bauthor{\bsnm{Cai},~\bfnm{Tony}\binits{T.}},
  \bauthor{\bsnm{Liu},~\bfnm{Weidong}\binits{W.}} \AND
  \bauthor{\bsnm{Luo},~\bfnm{Xi}\binits{X.}}
(\byear{2011}).
\btitle{A constrained {$\ell\sb 1$} minimization approach to sparse precision
  matrix estimation}.
\bjournal{J. Amer. Statist. Assoc.}
\bvolume{106}
\bpages{594--607}.
\bid{doi={10.1198/jasa.2011.tm10155}, issn={0162-1459}, mr={2847973}}
\bptok{imsref}%
\end{barticle}
\endbibitem

\bibitem[\protect\citeauthoryear{Cai, Liu and Zhou}{2012}]{CaiLiuZho}
\begin{bmisc}[auto:STB|2012/05/30|10:51:56]
\bauthor{\bsnm{Cai},~\bfnm{T.~T.}\binits{T.~T.}},
  \bauthor{\bsnm{Liu},~\bfnm{W.}\binits{W.}} \AND
  \bauthor{\bsnm{Zhou},~\bfnm{H.~H.}\binits{H.~H.}}
(\byear{2012}).
\bhowpublished{Optimal estimation of large sparse precision matrices.
  Unpublished manuscript}.
\bptok{imsref}%
\end{bmisc}
\endbibitem

\bibitem[\protect\citeauthoryear{Cand{\`e}s and Recht}{2009}]{CanRec09}
\begin{barticle}[mr]
\bauthor{\bsnm{Cand{\`e}s},~\bfnm{Emmanuel~J.}\binits{E.~J.}} \AND
  \bauthor{\bsnm{Recht},~\bfnm{Benjamin}\binits{B.}}
(\byear{2009}).
\btitle{Exact matrix completion via convex optimization}.
\bjournal{Found. Comput. Math.}
\bvolume{9}
\bpages{717--772}.
\bid{doi={10.1007/s10208-009-9045-5}, issn={1615-3375}, mr={2565240}}
\bptok{imsref}%
\end{barticle}
\endbibitem

\bibitem[\protect\citeauthoryear{Meinshausen and B{\"u}hlmann}{2006}]{MeiBuh06}
\begin{barticle}[mr]
\bauthor{\bsnm{Meinshausen},~\bfnm{Nicolai}\binits{N.}} \AND
  \bauthor{\bsnm{B{\"u}hlmann},~\bfnm{Peter}\binits{P.}}
(\byear{2006}).
\btitle{High-dimensional graphs and variable selection with the lasso}.
\bjournal{Ann. Statist.}
\bvolume{34}
\bpages{1436--1462}.
\bid{doi={10.1214/009053606000000281}, issn={0090-5364}, mr={2278363}}
\bptok{imsref}%
\end{barticle}
\endbibitem

\bibitem[\protect\citeauthoryear{Ravikumar et~al.}{2011}]{Ravetal11}
\begin{barticle}[mr]
\bauthor{\bsnm{Ravikumar},~\bfnm{Pradeep}\binits{P.}},
  \bauthor{\bsnm{Wainwright},~\bfnm{Martin~J.}\binits{M.~J.}},
  \bauthor{\bsnm{Raskutti},~\bfnm{Garvesh}\binits{G.}} \AND
  \bauthor{\bsnm{Yu},~\bfnm{Bin}\binits{B.}}
(\byear{2011}).
\btitle{High-dimensional covariance estimation by minimizing {$\ell\sb
  1$}-penalized log-determinant divergence}.
\bjournal{Electron. J. Stat.}
\bvolume{5}
\bpages{935--980}.
\bid{doi={10.1214/11-EJS631}, issn={1935-7524}, mr={2836766}}
\bptnote{check year}%
\bptok{imsref}%
\end{barticle}
\endbibitem

\end{thebibliography}
\end{document}